\documentclass[12pt]{amsart}

 \usepackage{amsfonts,graphics,amsmath,amsthm,amsfonts,amscd,amssymb,amsmath,latexsym,multicol}
\usepackage{epsfig,url}
\usepackage{flafter}

\makeatletter

\def\jobis#1{FF\fi
  \def\predicate{#1}%
  \edef\predicate{\expandafter\strip@prefix\meaning\predicate}%
  \edef\job{\jobname}%
  \ifx\job\predicate
}

\makeatother

\if\jobis{proposal}%
\else
\fi

 \usepackage[matrix, arrow]{xy}
 
\DeclareMathOperator{\Bs}{Bs}
\DeclareMathOperator{\Supp}{Supp}

\DeclareMathOperator{\codim}{codim}

\DeclareMathOperator{\Proj}{Proj}

\DeclareMathOperator{\exc}{exc}

\DeclareMathOperator{\Pic}{Pic}

 \newcommand{\C}{\mathbb C}
 \newcommand{\N}{\mathbb N}
 \newcommand{\PP}{\mathbb P}
 \newcommand{\Q}{\mathbb Q}
 \newcommand{\R}{\mathbb R}
 \newcommand{\Z}{\mathbb Z}
 \newcommand{\bir}{\dashrightarrow}

 \numberwithin{equation}{subsection}
 \numberwithin{footnote}{subsection}

 \newtheorem{cor}[subsubsection]{Corollary}
 \newtheorem{lem}[subsubsection]{Lemma}
 
 \newtheorem{thm}[subsubsection]{Theorem}
 \newtheorem{conj}[subsubsection]{Conjecture}

\newtheorem{prob}[subsubsection]{Problem}

{
 \newtheorem{defn}[subsubsection]{Definition}
 \newtheorem{exa}[subsubsection]{Example}
 \newtheorem{exa-cr}[subsection]{Example--Construction}

 \newtheorem{rem}[subsubsection]{Remark}

 \newtheorem{exe}[subsubsection]{Exercise}
}

 \newcommand{\ke}[1]{$\acute{\mbox{e}}$}
 \newcommand{\ku}[1]{$\acute{\mbox{u}$}}
 \newcommand{\kl}[1]{$\acute{\mbox{l}}$}
 \newcommand{\kh}[1]{$\acute{\mbox{h}}$}
 \newcommand{\kr}[1]{$\acute{\mbox{r}}$}
 \newcommand{\kx}[1]{$\acute{\mbox{x}}$}
 \newcommand{\ki}[1]{${\^\i}$}
 
 \newcommand{\toover}{\xrightarrow}

\title{Birational geometry}\thanks{3 April 2007. These are lecture notes for the \emph{Birational Geometry} course
which I taught at the University of Cambridge, Winter 2007.}
\author{Caucher Birkar}
\begin{document}
\maketitle

\tableofcontents
\vspace{2cm}
\clearpage

~
\vspace{2cm}
\section{\textbf{Introduction: overview}}
\vspace{1cm}

All varieties in this lecture are assumed to be algebraic over $\mathbb{C}$.\\

\emph{What classification means?} In every branch of mathematics, the problem of
classifying the objects arises naturally as the ultimate understanding of the
subject. If we have finitely many objects, then classification means trying to
see which object is isomorphic to another, whatever the isomorphism means. In this case,
the problem shouldn't be too difficult. However, if we have infinitely many objects, then
classification has a bit different meaning. Since each human can live for a finite
length of time, we may simply not have enough time to check every object in the theory.
On the other hand, objects in a mathematical theory are closely related, so one can
hope to describe certain properties of all the objects using only finitely many or
some of the objects which are nice in some sense.

For example let $V$ be a vector space over a field $k$. If $\dim V<\infty$, then
we can take a basis $\{v_1,\dots,v_n\}$ for this vector space and then every $v\in V$
can be uniquely written as $v=\sum a_iv_i$ where $a_i\in k$. We can say that we have
reduced the classification of elements of $V$ to the classification of the basis.
Formally speaking, this is what we want to do in any branch of mathematics but the
notion of "basis" and the basis "generating" all the objects could be very different.

In algebraic geometry, objects are varieties and relations are described using
maps and morphisms. One of the main driving forces of algebraic geometry is the following

\begin{prob}[Classification]
Classify varieties up to isomorphism.
\end{prob}

Some of the standard techniques to solve this problem, among other things, include
\begin{itemize}
\item Defining \emph{invariants} (e.g. genus, differentials, cohomology) so that
one can distinguish nonisomorphic varieties,
\item \emph{Moduli} techniques, that is, parametrising varieties by objects which are themselves
varieties,
\item \emph{Modifying} the variety (e.g. make it smooth) using certain operations
(eg birational, finite covers).
\end{itemize}

It is understood that this problem is too difficult even for curves so one needs to
try to solve a weaker problem.

\begin{prob}[Birational classification] Classify projective varieties up to birational isomorphism.
\end{prob}

By Hironaka's resolution theorem, each projective variety is \emph{birational} to a
smooth projective variety. So, we can try to classify
smooth projective varieties up to birational isomorphism. In fact, smooth birational
varieties have good common properties such as common plurigenera,
Kodaira dimension and irregularity.

In this stage, one difficulty is that in each birational class, except in dimension one,
there are too many smooth varieties. So, one needs to look for a more special and
subtle representative.

\begin{exa}[Curves] Projective curves are one dimensional projective varieties
(i.e. compact Riemann surfaces). Curves $X$ in $\mathbb{P}^2$ are defined by a single homogeneous
polynomial $F$.
A natural and important "invariant" is the \emph{degree} defined as $\deg(X)=\deg F$. The degree
is actually not an invariant because a line and a conic have different degrees but they
could be isomorphic.
However, using the degree we can simply define an invariant: \emph{genus}, which is defined as
$$
g(X)=\frac{1}{2}(\deg(X)-1)(\deg(X)-2)
$$
Topologically, $g(X)$ is the number of handles on $X$.\\

 \emph{Moduli spaces.} For each $g$, there is a moduli space $\mathcal{M}_g$ of
smooth projective curves of genus $g$.

\begin{displaymath}
g(X) = \left\{ \begin{array}{ll}
0 & \textrm{iff $X\simeq \mathbb{P}^1$}\\
1 & \textrm{iff $X$ is elliptic}\\
\ge 2 & \textrm{iff $X$ is of general type}
\end{array} \right.
\end{displaymath}

 These correspond to
 \begin{displaymath}
g(X) = \left\{ \begin{array}{ll}
0 & \textrm{iff $X$ has positive curvature}\\
1 & \textrm{iff $X$ has zero curvature}\\
\ge 2 & \textrm{iff $X$ has negative curvature}
\end{array} \right.
\end{displaymath}

On the other hand, $g(X)=h^0(X,\omega_X)$ where $\omega_X=\mathcal{O}_X(K_X)$ is the canonical sheaf,
that is, $\Omega^d$.

 \begin{displaymath}
g(X) = \left\{ \begin{array}{ll}
0 & \textrm{iff $\deg K_X<0$}\\
\ge 1 & \textrm{iff $\deg K_X\ge 0$}
\end{array} \right.
\end{displaymath}

\end{exa}

\begin{defn} For a smooth projective variety $X$, define the $m$-th plurigenus as
$$
P_m(X):=h^0(X,\omega_X^{\otimes m})
$$
Note that $P_1(X)=g(X)$.
Define the Kodaira dimension of $X$ as
$$
\kappa(X):=\limsup_{m\to\infty} \frac{\log P_m(X)}{\log m}
$$
\end{defn}

If $\dim X=d$, then $\kappa(X)\in\{-\infty, 0, 1, \dots, d\}$.
Moreover, the Kodaira dimension and the plurigenera $P_m(X)$ are all birational invariants.\\

\begin{exa} If $\dim X=1$, then
 \begin{displaymath}
\kappa(X) = \left\{ \begin{array}{ll}
-\infty & \textrm{iff $\deg K_X<0$}\\
\ge 0 & \textrm{iff $\deg K_X\ge 0$}
\end{array} \right.
\end{displaymath}
\vspace{0.1cm}
\end{exa}

In higher dimension, the analogue of canonical with negative degree is the notion of
   a \emph{Mori fibre space} defined as a fibre type contraction $Y\to Z$ which is a
 $K_Y$-negative extremal fibration with connected fibres.
But the analogue of canonical divisor with nonnegative degree is the notion of a \emph{minimal}
variety defined as $Y$ having $K_Y$ nef, i.e., $K_Y\cdot C\ge 0$ for any curve $C$ on $Y$.

\begin{conj}[Minimal model] Let $X$ be a smooth projective variety.\\\\ Then,
\begin{itemize}
\item If $\kappa(X)=-\infty$, then $X$ is birational to a Mori fibre space $Y\to Z$.\\
\item If $\kappa(X)\ge 0$, then $X$ is birational to a minimal variety $Y$.
\end{itemize}
\end{conj}

\begin{conj}[Abundance] Let $Y$ be a minimal variety.\\\\ Then,
there is an \emph{Iitaka fibration} $\phi\colon Y\to S$ with connected fibres
 and an ample divisor $H$ on $S$ such that
$$
K_Y=\phi^* H
$$

 In fact, \\

$\bullet$ $\phi(C)=\rm pt. \Longleftrightarrow  K_Y\cdot C=0$ for any curve $C$ on $Y$.\\

$\bullet$ $\dim S=\kappa(Y)$.
\end{conj}

\begin{exa}[Enriques-Kodaira classification of surfaces]
 The minimal model conjecture and the abundance conjecture hold for surfaces. More precisely,\\
$\kappa(X)=-\infty \implies$  $X$ is birational to $\PP^2$ or a $\PP^1$-bundle over some curve.\\
$\kappa(X)=0 \implies$ $X$ is birational to a K3 surface, an Enriques surface or an \'etale
quotient of an abelian surface. The Iitaka fibration $\phi\colon Y\to S=\rm pt.$ is trivial.\\
$\kappa(X)=1 \implies$ $X$ is birational to a minimal elliptic surface $Y$. The Iitaka fibration
$\phi\colon Y\to S$ is an elliptic fibration.\\
$\kappa(X)=2 \implies$ $X$ is birational to a minimal model $Y$. The Iitaka fibration
$\phi\colon Y\to S$ is birational.
\end{exa}

\begin{rem}[Classical MMP]\label{rem-classical MMP} To get the above classification
one can use the classical \emph{minimal model program} (MMP) as
follows. If there is a $-1$-curve $E$ (i.e. $E\simeq \PP^1$ and $E^2=-1$) on $X$,
then by Castelnuovo theorem we can contract $E$ by a birational morphism $f\colon X\to X_1$
where $X_1$ is also smooth. Now replace $X$ with $X_1$ and continue the process. In each step,
the Picard number $\rho(X)$ drops by $1$. Therefore, after finitely many steps, we get a
smooth projective variety $Y$ with no $-1$-curves. Such $Y$ turns out to be among one of the
classes in Enriques classification.
\end{rem}

In higher dimension (ie. $\dim \ge 3$), we would like to have a program similar to the
classical MMP for surfaces. However, this is not without difficulty. Despite being a
very active area of algebraic geometry and large number of people working on this
program for nearly three decades, there are still many fundamental problems to be solved.
For the precise meaning of the terminology see other sections of the lecture. In any case,
some of the main problems including the minimal model conjecture and abundance
conjecture are:
 
\begin{itemize}
\item \emph{Minimal model conjecture.} Open in dimension $\ge 5$.
\item \emph{Abundance conjecture.} Open in dimension $\ge 4$.
\item \emph{Flip conjecture:} Flips exist. Recently solved.
\item \emph{Termination conjecture:} Any sequence of flips terminates. Open in dimension $\ge 4$.
\item \emph{Finite generation conjecture:} For any smooth projective variety $X$, the canonical ring
$$
R(X):=\bigoplus_{m=0}^{\infty} H^0(X,\omega_X^{\otimes m})
$$
is finitely generated. Recently solved.
\item \emph{Alexeev-Borisov conjecture:} Fano varieties of dimension $d$ with bounded singularities, are bounded.
Open in dimension $\ge 3$.
\item \emph{ACC conjectures on singularities:} Certain invariants of singularities (i.e. minimal log discrepancies and
log canonical thresholds) satisfy the ascending chain condition (ACC). Open in dimension $\ge 3$ and $\ge 4$
respectively.
\item \emph{Uniruledness:} Negative Kodaira dimension is equivalent to uniruledness. Open in dimension $\ge 4$.
\end{itemize}

The following two theorems and their generalisations are the building blocks of the
techniques used in minimal model program.

\begin{thm}[Adjunction]
Let $X$ be a smooth variety and $Y\subset X$ a smooth subvariety.
Then,
$$
(K_X+Y)\vert_Y=K_Y
$$
\end{thm}

\begin{thm}[Kodaira vanishing]
Let $X$ be a smooth projective variety and $H$ an ample divisor on $X$. Then,
$$
H^i(X,K_X+H)=0
$$
for any $i>0$.
\end{thm}

\clearpage

~ \vspace{2cm}
\section{\textbf{Preliminaries}}
\vspace{1cm}

\begin{defn}
In this lecture, by a \emph{variety} we mean an irreducible quasi-projective variety over $\C$.
Two varieties $X,X'$ are called \emph{birational} if there is a rational map $f\colon X\bir X'$
which has an inverse, i.e., if $X,X'$ have isomorphic open subsets.
\end{defn}

\begin{defn}[Contraction]
A contraction is a projective morphism $f\colon X\to Y$ such that
$f_*\mathcal{O}_X=\mathcal{O}_Y$ (so, it has connected fibres).
\end{defn}

\begin{exa}[Zariski's main theorem]
Let $f\colon X\to Y$ be a projective birational morphism where $Y$ is normal.
Then, $f$ is a contraction.
\end{exa}

\begin{exa} Give a finite map which is not a contraction.
\end{exa}

\begin{rem}[Stein factorisation]
Let $f\colon X\to Y$ be a projective morphism. Then, it can be factored
through $g\colon X\to Z$ and $h\colon Z\to Y$ such that $g$ is a contraction
and $h$ is finite.
\end{rem}

\begin{defn}
Let $X$ be a normal variety. A \emph{divisor} (resp. \emph{$\Q$-divisor, $\R$-divisor}) is
$\sum_{i=1}^md_iD_i$ where $D_i$ are prime divisors and $d_i\in\Z$ (resp.
$d_i\in\Q$, $d_i\in\R$). A $\Q$-divisor $D$ is called \emph{$\Q$-Cartier} if
$mD$ is Cartier for some $m\in\N$. Two $\Q$-divisors $D,D'$ are called
\emph{$\Q$-linearly equivalent}, denoted by $D\sim_{\Q}D'$, if $mD\sim mD'$
for some $m\in\N$. $X$ is called \emph{$\Q$-factorial} if
every $\Q$-divisor is $\Q$-Cartier.

\end{defn}

\begin{defn}
Let $X$ be a normal variety and $D$ a divisor on $X$. Let $U=X-X_{sing}$ be the smooth
subset of $X$ and $i\colon U\to X$ the inclusion. Since $X$ is normal,
$\dim X_{sing}\le \dim X-2$. So, every divisor on $X$ is uniquely determined
by its restriction to $U$. In particular, we define the \emph{canonical divisor} $K_X$
of $X$ to be the closure of the canonical divisor $K_U$.

For a divisor $D$, we associate the sheaf $\mathcal{O}_X(D):=i_*\mathcal{O}_U(D)$.
It is well-known that $\mathcal{O}_X(K_X)$ is the same as the dualising sheaf $\omega_X$
in the sense of [\ref{Hart}, III.7.2] if $X$ is projective. See [\ref{KM}, Proposition 5.75].
\end{defn}

\begin{exa}
The canonical divisor of $\PP^n$ is just $-(n+1)H$ where $H$ is a hyperplane.
\end{exa}

\begin{defn}
Let $X$ be a projective variety. A $\Q$-Cartier divisor $D$ on $X$ is called
\emph{nef} if $D\cdot C\ge 0$ for any curve $C$ on $X$. Two $\Q$-Cartier divisors $D,D'$ are called
\emph{numerically equivalent}, denoted by $D\equiv D'$, if $D\cdot C=D'\cdot C$ for
any curve $C$ on $X$. Now let $V,V'$ be two $\R$-1-cycles on $X$. We call them numerically
equivalent, denoted by $V\equiv V'$ if $D\cdot V=D\cdot V'$ for any $\Q$-Cartier divisor
$D$ on $X$.
\end{defn}

\begin{defn}
Let $X$ be a projective variety. We define
 $$N_1(X):=\mbox{group of $\R$-1-cycles}/\equiv$$
 $$NE(X):=\mbox{the cone in $N_1(X)$ generated by effective $\R$-1-cycles}$$
 $$\overline{NE}(X):=\mbox{closure of $NE(X)$ inside $N_1(X)$} $$
 $$N^1(X):=(\Pic(X)/\equiv)\otimes_{\Z} \R$$

It is well-known that $\rho(X):=\dim N_1(X)=\dim N^1(X)<\infty$ which is called
the \emph{Picard number} of $X$.
\end{defn}

\begin{defn}
Let $C\subset \R^n$ be a cone with the vertex at the origin.
A subcone $F\subseteq C$ is called an \emph{extremal face} of $C$ if
for any $x,y\in C$, $x+y\in F$ implies that $x,y\in F$. If $\dim F=1$,
we call it an \emph{extremal ray}.
\end{defn}

\begin{thm}[Kleiman ampleness criterion]
Let $X$ be a projective variety and $D$ a $\Q$-Cartier divisor. Then,
$D$ is ample iff $D$ is positive on $\overline{NE}(X)-\{0\}$.
\end{thm}

\begin{thm}
Let $X$ be a projective surface and $C$ an irreducible curve. (i) If $C^2\le 0$,
then $C$ is in the boundary of $\overline{NE}(X)$. (ii) Moreover, if $C^2<0$, then
$C$ generates an extremal ray.
\end{thm}
\begin{proof}
 See [\ref{KM}, Lemma 1.22].
\end{proof}

\begin{exa}
Let $X$ be a smooth projective curve. Then we have a natural exact sequence
$0\to \Pic^0(X) \to \Pic(X) \to \Z\to 0$. So, $\Pic/\Pic^0\simeq \Z$. Therefore,
$N_1(X)\simeq \R$ and $\overline{NE}(X)$ is just $\R_{\ge 0}$.
If $X=\PP^n$, we have a similar story.
\end{exa}

\begin{exa}
Let $X=\PP^1\times \PP^1$ be the quadric surface. $N_1(X)\simeq \R^2$ and $\overline{NE}(X)$
has two extremal rays, each one is generated by fibres of one of the two natural projections.

If $X$ is the blow up of
$\PP^2$ at a point $P$, then $N_1(X)\simeq \R^2$ and $\overline{NE}(X)$ has two extremal
rays. One is generated by the exceptional curve of the blow up and the other one
is generated by the birational transform (see below) of all the lines passing through $P$.

If $X$ is a cubic surface, it is well-known that it contains exactly 27 lines,
$N_1(X)\simeq \R^7$ and $\overline{NE}(X)$ has 27 extremal rays, each one generated by
one of the 27 lines.
 
If $X$ is an abelian surface, one can prove that $\overline{NE}(X)$ has a round shape,
that is it does not look like polyhedral. This happens because an abelian variety
is in some sense homogeneous.

Finally, there are surfaces $X$ which have infinitely many $-1$-curves. So,
they have infinitely many extreml rays. To get such a surface one can blow up the
projective plane at nine points which are the base points of a general pencil of cubics.
\end{exa}

We will mostly be interested in the extremal rays on which $K_X$ is negative.

\begin{defn}[Exceptional set]
Let $f\colon X\to Y$ be a birational morphism of varieties. $\exc(f)$ is the set of those
$x\in X$ such that $f^{-1}$ is not regular at $f(x)$.
\end{defn}

\begin{defn}
Let $f\colon X\bir Y$ be a birational map of normal varieties and $V$ a prime cycle on $X$.
Let $U\subset X$ be the open subset where $f$ is regular. If $V\cap U\neq \emptyset$, then define
the \emph{birational transform} of $V$ to be the closure of $f(U\cap V)$ in $Y$. If
$V=\sum a_iV_i$ is a cycle and $U\cap V_i\neq \emptyset$, then the birational transform
of $V$ is defined to be $\sum a_iV_i^{\sim}$ where $V_i^{\sim}$ is the birational transform
of the prime component $V_i$.
\end{defn}

\begin{defn}[Log resolution]
Let $X$ be a variety and $D$ a $\Q$-divisor on $X$. A projective birational morphism
$f\colon Y\to X$ is a log resolution of $X,D$ if $Y$ is smooth,
$\exc(f)$ is a divisor and $\exc(f)\cup f^{-1}(\Supp D)$ is a simple normal crossing divisor.
\end{defn}

\begin{thm}[Hironaka]
Let $X$ be a variety and $D$ a $\Q$-divisor on $X$. Then, a log resolution exists.
\end{thm}

\clearpage

~ \vspace{2cm}
\section{\textbf{Singularities}}
\vspace{1cm}

\begin{defn}
Let $X$ be a variety and $D$ a divisor on $X$. We call $X,D$ \emph{log smooth} if
$X$ is smooth and the components of $D$ have simple normal crossings.
A \emph{pair} $(X,B)$ consists of a normal variety $X$ and a $\Q$-divisor $B$ with
coefficients in $[0,1]$ such that $K_X+B$ is $\Q$-Cartier.
Now let $f\colon Y\to X$ be a log resolution of a pair $(X,B)$. Then, we can write

$$
K_Y=f^*(K_X+B)+A
$$

 For a prime divisor $E$ on $Y$, we define the \emph{discrepancy} of $E$ with respect
 to $(X,B)$ denoted by $d(E,X,B)$ to be the coefficient of $E$ in $A$. Note that
if $E$ appears as a divisor on any other resolution, then $d(E,X,B)$ is the same.
\end{defn}

\begin{rem}[Why pairs?]
The main reason for considering pairs is the various kinds of adjunction, that is,
relating the canonical divisor of two varieties which are closely related. We have already
seen the adjunction formula $(K_X+B)\vert_S=K_S$ where $X,S$ are smooth and $S$ is a prime
divisor on $X$. It is natural to consider $(X,S)$ rather than just $X$.

Now let $f\colon X\to Z$ be a finite morphism. It often happens that $K_X=f^*(K_Z+B)$
for some boundary $B$. Similarly, when $f$ is a fibration and $K_X\sim_{\Q}f^*D$ for
some $\Q$-Cartier divisor $D$ on $Z$, then under good conditions $K_X\sim_{\Q}f^*(K_Z+B)$
for some boundary $B$ on $Z$. Kodaira's canonical bundle formula for
an elliptic fibration of a surface is a clear example.
\end{rem}

\begin{rem}
Let $X$ be a smooth variety and $D$ a $\Q$-divisor on $X$. Let $V$ be a
smooth subvariety of $X$ of codimension $\ge 2$ and $f\colon Y\to X$ the
blow up of $X$ at $V$, and $E$ the exceptional divisor. Then, the coefficient
of $E$ in $A$ is $\codim V -1-\mu_V D$ where

$$
K_Y=f^*(K_X+D)+A
$$
and $\mu$ stands for multiplicity.
\end{rem}

\begin{defn}[Singularities]
Let $(X,B)$ be a pair. We call it \emph{terminal} (resp. \emph{canonical}) if $B=0$ and
there is a log resolution $f\colon Y\to X$ for which $d(E,X,B)>0$ (resp. $\ge 0$) for
any exceptional prime divisor $E$ of $f$. We call the pair \emph{Kawamata log terminal}
 ( resp. \emph{log canonical}) if there is a log resolution $f$ for which $d(E,X,B)>-1$ (resp. $\ge -1$)
 for any prime divisor $E$ on $Y$ which is exceptional for $f$ or the birational transform of
 a component of $B$. The pair is called \emph{divisorially log terminal} if there is a log resolution
 $f$ for which $d(E,X,B)>-1$ for any exceptional prime divisor $E$ of $f$. We
 usually use abbreviations klt, dlt and lc for Kawamata log terminal, divisorially log terminal
 and log canonical respectively.
\end{defn}

\begin{lem}
Definition of all kind of singularities (except dlt) is independent of the choice of the log resolution.
\end{lem}
\begin{proof}
Suppose that $(X,B)$ is lc with respect to a log resolution $(Y,B_Y)$. Let $(Y',B_Y')$ be
another log resolution and $(Y'',B_Y'')$ a common blowup of $Y$ and $Y'$. The remark above
shows that in this way we never get a discrepancy less than $-1$.

Other cases are similar.
\end{proof}

\begin{exe} Prove that if $(X,B)$ is not lc and $\dim X>1$, then for any integer $l$ there is $E$
such that $d(E,X,B)<l$.
\end{exe}

\begin{exe} Prove that a smooth variety is terminal.
\end{exe}

\begin{exe} Prove that terminal $\implies$ canonical $\implies$ klt $\implies$ dlt $\implies$ lc.
\end{exe}

\begin{exe} If $(X,B+B')$ is terminal (resp. canonical, klt, dlt or lc) then so is $(X,B)$ where
$B'\ge 0$ is $\Q$-Cartier.
\end{exe}

\begin{exe}
 Let $(X,B)$ be a pair and $f\colon Y\to X$ a log resolution. Let $B_Y$ be the divisor on $Y$
 for which $K_Y+B_Y=f^*(K_X+B)$. Prove that $(X,B)$ is\\
$\bullet$ terminal iff $B_Y\le 0$ and $\Supp B_Y=\exc(f)$,\\
$\bullet$ canonical iff $B_Y\le 0$,\\
$\bullet$ klt iff each coefficient of $B_Y$ is $<1$,\\
$\bullet$ lc iff each coefficient of $B_Y$ is $\le 1$.
\end{exe}

\begin{exa}
Let $(X,B)$ be a pair of dimension $1$. Then, it is lc (or dlt) iff each coefficient of $B$ is
$\le 1$. It is klt iff each of coefficient of $B$ is $<1$. It is canonical (or terminal) iff
$B=0$.
\end{exa}

\begin{exa}
When $(X,B)$ is log smooth then we have the most simple kind of singularities. It is
easy to see what type of singularities this pair has.

Now let
$(\PP^2,B)$ be a pair where $B$ is a nodal curve. This pair is lc but not dlt.
However, the pair $(\PP^2,B)$ where $B$ is a cusp curve is not lc.
\end{exa}
 
\begin{exa}
  Let's see what terminal, etc. mean for some of the simplest
surface singularities. Let $X$ be a smooth surface containing
a curve $E=\PP^1$ with $E^2 = -a$, $a>0$. (Equivalently, the normal
bundle to $E$ in $X$ has degree $-a$.) It's known that one can
contract $E$ to get a singular surface $Y$, $f\colon X\to Y$. Explicitly,
$Y$ is locally analytically the cone over the rational normal curve
    $\PP^1 < \PP^a$;
so for $a=1$, $Y$ is smooth, and for $a=2$, $Y$ is the surface node
${x^2+y^2-z^2=0}$ in $\mathbb{A}^3$ (a canonical singularity).

Then $K_E = (K_X+E)\vert_E$, and $K_E$ has degree $-2$ on $E$
since $E$ is isomorphic to $\PP^1$, so
   $K_X \cdot E = -2+a$.
This determines the discrepancy c in
   $K_X = f^*(K_Y) + cE$,
because $f^*(K_Y)\cdot E = K_Y\cdot (f_*(E)) = K_Y \cdot 0=0$. Namely,
  $c = (a-2)/(-a)$.
So for $a=1, c=1$ and $Y$ is terminal (of course, since it's smooth);
for $a=2, c=0$ and $Y$ is canonical (here $Y$ is the node); and for $a\ge 3$,
$c$ is in $(-1,0)$ and $Y$ is klt.

Just for comparison:
if you contract a curve of genus 1, you get an lc singularity
which is not klt; and if you contract a curve of genus at least 2,
it is not even lc.
\end{exa}

\begin{defn}
Let $f\colon X\to Z$ be a projective morphism of varieties
and $D$ a $\Q$-Cartier divisor on $X$. $D$ is called \emph{nef over $Z$} if
$D\cdot C\ge 0$ for any curve $C\subseteq X$ contracted by $f$.
$D$ is called \emph{numerically zero over $Z$} if $D\cdot C\ge 0$ for any curve
$C\subseteq X$ contracted by $f$.
\end{defn}

\begin{lem}[Negativity lemma]
Let $f\colon Y\to X$ be a projective birational morphism of normal varieties.
Let $D$ be a $\Q$-Cartier divisor on $Y$ such that $-D$ is nef over $X$. Then,
$D$ is effective iff $f_*D$ is.
\end{lem}
\begin{proof}
First by localising the problem and taking hyperplane sections we can assume
that $X,Y$ are surfaces and $D$ is contracted by $f$ to a point $P\in X$.
Now by replacing $Y$ with a
resolution and by taking a Cartier divisor $H$ passing through $P$, we can find
an effective exceptional divisor $E$ which is antinef/$X$.

Let $e$ be the minimal non-negative number for which $D+eE\ge 0$. If $D$ is not
effective, then $D+eE$ has coefficient zero at some component. On the other hand,
$E$ is connected. This is a contradiction.
\end{proof}

\begin{defn}[Minimal resolution]
Let $X$ be a normal surface and $f\colon Y\to X$ a resolution. $f$ or $Y$ is called
the \emph{minimal resolution} of $X$ if $f'\colon Y'\to X$ is any other resolution,
then $f'$ is factored through $f$.
\end{defn}

\begin{rem} The matrix $[E_i\cdot E_j]$ is negative definite  for the resolution of a normal
surface singularity. So, we can compute discrepancies by intersecting each exceptional curve with
$K_Y+B_Y=f^*K_X$.
\end{rem}

\begin{lem}
Let $f\colon Y\to X$ be the minimal resolution of a normal surface $X$. Then,
$K_Y+B_Y=f^*K_X$ where $B_Y\ge 0$.
\end{lem}
\begin{proof}
First prove that $K_Y$ is nef over $X$ using the formula $(K_Y+C)\cdot C=2p_a(C)-2$
for a proper curve $C$ on $Y$. Then, the negativity lemma implies that $B_Y\ge 0$.
\end{proof}

\begin{thm}
A surface $X$ is terminal iff it is smooth.
\end{thm}
\begin{proof}
If $X$ is smooth then it is terminal. Now suppose that $X$ is terminal
and let $Y\to X$ be a minimal resolution and let $K_Y+B_Y$ be the pullback of
$K_X$. Since $X$ is terminal, $d(E,X,0)>0$ for any exceptional divisor. Thus,
$B_Y<0$ a contradiction. So, $Y=X$.
\end{proof}

\begin{cor}\label{cor-terminal} By taking hyperplane sections, one can show that terminal
varieties are smooth in codimension two.
\end{cor}

\begin{rem}
The following are equivalent:\\
 $X$ has canonical surface singularities (:= Du Val singularities)\\
 Locally analytically $X$ is given by the following equations in $\mathbb{A}^3$:\\
 A: $x^2+y^2+z^{n+1}=0$\\
 D: $x^2+y^2z+z^{n-1}=0$\\
 E$_6$: $x^2+y^3+z^4=0$\\
 E$_7$: $x^2+y^3+yz^3=0$\\
 E$_8$: $x^2+y^3+z^5=0$\\
\end{rem}

\begin{lem}
If $X$ is klt, then all exceptional curves of the minimal resolution are smooth rational
curves.
\end{lem}
\begin{proof}
Let $E$ be a an exceptional curve appearing on the minimal resolution $Y$. Then, $(K_Y+eE)\cdot E\le 0$
for some $e<1$. So, $(K_Y+E+(e-1)E)\cdot E=2p_a(E)-2+(e-1)E^2\le 0$ which in turn implies that
$p_a(E)\le 0$. Therefore, $E$ is a smooth rational curve.
\end{proof}

\begin{rem} the dual graph of a surface klt singularity: A, D and E type [\ref{Prokhorov}].
\end{rem}

\begin{exa} \emph{Singularities in higher dimension}.
Let $X$ be defined by $x^2+y^2+z^2+u^2=0$ in $\mathbb{A}^4$. Then, by blowing up the
origin of $\mathbb{A}^4$ we get a resolution $Y\to X$ such that we have a
single exceptional divisor $E$ isomorphic to the quadric surface $\PP^1\times \PP^1$.
Suppose that $K_Y=f^*K_X+eE$. Take a fibre $C$ of the projection $E\to \PP^1$.
Either by calculation or more advanced methods, one can show that $K_Y\cdot C<0$ and
$E\cdot C<0$. Therefore, $e>0$. So, it is a terminal singularity.
\end{exa}

\begin{rem}[Toric varieties]  Suppose that $X$ is the toric variety
associated to a cone $\sigma\subset N_{\R}$,
\begin{itemize}
\item $X$ is smooth iff $\sigma$ is regular, that is primitive generators of each face of $\sigma$ consists
of a part of a basis of $N$,
\item $X$ is $\Q$-factorial iff $\sigma$ is simplicial,
\item $X$ is terminal iff $\sigma$ is terminal, that is, there is $m\in M_{\Q}$ such that
$m(P)=1$ for each primitive generator $P\in \sigma \cap N$, and $m(P)>1$ for any other
$P\in N\cap \sigma -\{0\}$,
\item If $K_X$ is $\Q$-Cartier, then $X$ is klt.
\end{itemize}

See [\ref{Dais}] and [\ref{Matsuki}] for more information.
\end{rem}

\clearpage
~ \vspace{2cm}
\section{\textbf{Minimal model program}}
\vspace{1cm}

\subsection{Kodaira dimension}

\begin{defn}\label{def-D-morphism}
Let $D$ be a divisor on a normal projective variety $X$. If $h^0(D)=n\neq 0$, then we define
the rational map $\phi_{D}\colon X \bir \PP^{(n-1)}$ as
$$
\phi_{D}(x)=(f_1(x):\dots:f_n(x))
$$
where $\{f_1,\dots,f_n\}$ is a basis for $H^0(X,D)$.
\end{defn}
 
\begin{defn}[Kodaira dimension] For a divisor $D$ on a normal projective variety $X$,
define
$$
\kappa(D)=\max\{\dim \phi_{mD}(X)\mid h^0(X,mD)\neq 0\}
$$
if $h^0(X,mD)\neq 0$ for some $m\in \N$, and $-\infty$ otherwise.

If $D$ is a $\Q$-divisor, then define $\kappa(D)$ to be $\kappa(lD)$
where $lD$ is integral. In particular, by definition $\kappa(D)\in\{-\infty, 0,\dots,\dim X\}$.

For a pair $(X,B)$, define the Kodaira dimension $\kappa(X,B):=\kappa(K_X+B)$.
\end{defn}

 \begin{defn}[Big divisor]
A $\Q$-divisor on a normal variety $X$ of dimension $d$ is \emph{big} if $\kappa(D)=d$. In particular,
if $\kappa(X,B)=d$, we call $(X,B)$ of general type.
\end{defn}

\begin{exe}\label{e-Kodaira}
Prove that the definition of the Kodaira dimension of a $\Q$-divisor is well-defined. That is,
$\kappa(D)=\kappa(lD)$ for any $l\in\N$ if $D$ is integral.
\end{exe}

\begin{exe}
Let $H$ be an ample divisor on a normal projective variety $X$. Prove that $\kappa(H)=\dim X$.
\end{exe}

\begin{exe}
Let $f\colon Y\to X$ be a contraction of normal projective varieties and $D$ a $\Q$-Cartier divisor on $X$.
Prove that $\kappa(D)=\kappa(f^*D)$.
\end{exe}

\begin{exe}
Let $D$ be a divisor on a normal projective variety $X$. Prove that,\\
$\bullet$ $\kappa(D)=-\infty$ $\Longleftrightarrow$ $h^0(mD)=0$ for any $m\in\N$.\\
$\bullet$ $\kappa(D)=0$  $\Longleftrightarrow$ $h^0(mD)\le 1$ for any $m\in\N$ with equality for some $m$.\\
$\bullet$ $\kappa(D)\ge 1$ $\Longleftrightarrow$ $h^0(mD)\ge 2$ for some $m\in\N$.

\end{exe}

\begin{exe}*
Let $D,L$ be $\Q$-divisors on a normal projective variety $X$ where $D$ is big. Prove that
there is $m\in\N$ such that $h^0(X,mD-L)\neq 0$.
\end{exe}

\begin{exe}* For a $\Q$-divisor $D$ on a projective normal variety $X$,
define
$$
P_m(D):=h^0(X,\lfloor mD\rfloor)
$$
Define the Kodaira dimension of $D$ as
$$
\kappa(D):=\limsup_{m\to\infty} \frac{\log P_m(D)}{\log m}
$$

Now prove that this definition is equivalent to the one given above.
\end{exe}

\subsection{Basics of minimal model program}

\begin{defn}[Minimal model-Mori fibre space]
A projective lc pair $(Y,B_Y)$ is called \emph{minimal} if $K_Y+B_Y$ is nef. A $(K_Y+B_Y)$-negative
extremal contraction $g\colon Y\to Z$ is called a \emph{Mori fibre space} if $\dim Y>\dim Z$.
 
 Let $(X,B)$, $(Y,B_Y)$ be lc pairs and $f\colon X\bir Y$ a birational map whose inverse does not
 contract any divisors such that $B_Y=f_*B$. $(Y,B)$ is called a \emph{minimal model} for $(X,B)$
 if $K_Y+B_Y$ is nef and if $d(E,X,B)>d(E,Y,B_Y)$ for any prime divisor $E$ on $X$ contracted by $f$.
 
 $g\colon Y\to Z$ is a \emph{Mori fibre space} for $(X,B)$ if it is a $(K_Y+B_Y)$-negative
 extremal contraction such that $d(E,X,B)>d(E,Y,B_Y)$ for any prime divisor $E$ on $X$ contracted by $f$.
\end{defn}

\begin{conj}[Minimal model] Let $(X,B)$ be a projective lc pair.
Then,\\
$\bullet$ If $\kappa(X,B)=-\infty$, then $(X,B)$ has a Mori fibre space.\\
$\bullet$ If $\kappa(X,B)\ge 0$, then $(X,B)$ has a minimal model.
\end{conj}

\begin{conj}[Abundance] Let $(Y,B)$ be a minimal lc pair.
Then, $m(K_Y+B)$ is base point free for some $m\in\N$. This, in particular,
means that there is a contraction $h\colon Y\to S$ called the \emph{Iitaka fibration}
and an ample $\Q$-divisor $H$ on $S$ such that
$$
K_Y+B=h^* H
$$

 In fact, \\

$\bullet$ $h(C)=\rm pt. \Longleftrightarrow  (K_Y+B)\cdot C=0$ for any curve $C$ on $Y$.\\

$\bullet$ $\dim S=\kappa(Y,B)$.
\end{conj}

\begin{conj}[Iitaka]
Let $f\colon X\to Z$ be a contraction of smooth projective varieties with smooth general fibre $F$.
Then,
$$
\kappa(X)\ge \kappa(Z)+\kappa(F)
$$
\end{conj}

\begin{defn}[Contraction of an extremal ray]
Let $R$ be an extremal ray of a normal projective variety $X$. A contraction $f\colon X\to Z$ is
the contraction of $R$ if
$$
f(C)=pt. \Longleftrightarrow [C]\in R
$$
for any curve $C\subset X$.
\end{defn}

\begin{rem}[Types of contraction]
For the contraction of an extremal ray $R$ we have the following possibilities:
\begin{description}
 \item[Divisorial] $f$ is birational and contracts divisors.
 \item[Flipping]   $f$ is birational and does not contract divisors.
 \item[Fibration]  $f$ is not birational.
\end{description}
\end{rem}

\begin{defn}[Flip]
 Let $(X,B)$ be a pair where $X$ is projective. A $(K_X+B)$-flip is a diagram

$$
\xymatrix{ X\ar[rd]_{f}& \dashrightarrow & X^+\ar[ld]^{f^+} \\
&Z&}
$$

such that

\begin{itemize}
\item $X^+$ and $Z$ are normal varieties.
\item $f$ and $f^+$ are small projective birational contractions, where small
means that they contract subvarieties of codimension $\geq 2$.
\item $f$ is the contraction of an extremal ray $R$.
\item $-(K_X+B)$ is ample over $Z$, and $K_{X^+}+B^+$ is ample over $Z$
where $B^+$ is the birational transform of $B$.
 \end{itemize}
\end{defn}

\begin{exe}
Suppose that $X$ is $\Q$-factorial and projective. Prove that $\Q$-factoriality is preserved after divisorial
contractions and flips.
\end{exe}

\begin{defn}[Minimal model program: MMP]
The minimal model program can be described in different level of generality. The
following seems to be reasonable.
Let $(X,B)$ be a dlt pair where $X$ is $\Q$-factorial and projective. The following process is called
the minimal model program if exists:

If $K_X+B$ is not nef, then there is an $(K_X+B)$-extremal ray $R$ and its contraction $f\colon X\to Z$.
If $\dim Z<\dim X$, then we get a Mori fibre space and we stop. If $f$ is a divisorial contraction,
we replace $(X,B)$ with $(Z,f_*B)$ and continue. If $f$ is flipping, we replace $(X,B)$ with the
flip $(X^+,B^+)$ and continue. After finitely many steps, we get a minimal model or a Mori fibre space.
\end{defn}

\begin{exa}
Classical MMP for smooth projective surfaces.
\end{exa}

\begin{conj}[Termination]
Let $(X,B)$ be a dlt pair where $X$ is $\Q$-factorial and projective. Any sequence of $(K_X+B)$-flips
terminates.
\end{conj}

\subsection{Cone and contraction, vanishing, nonvanishing, and base point freeness}

\begin{thm}[Cone and contraction]
Let $(X,B)$ be a klt pair where $X$ is projetive.
Then, there is a set of $(K_X+B)$-negative extremal rays $\{R_i\}$ such that\\
$\bullet$ $\overline{NE}(X)=\overline{NE}(X)_{K_X+B\ge 0}+\sum_iR_i$.\\
$\bullet$ Each $R_i$ contains the class of some curve $C_i$, that is $[C_i]\in R_i$.\\
$\bullet$ $R_i$ can be contracted.\\
$\bullet$ $\{R_i\}$ is discrete in $\overline{NE}(X)_{K_X+B< 0}$.
\end{thm}

\begin{rem} Remember that a divisor $D$ on a normal variety $X$ is called \emph{free} if its base locus
$$
\Bs |D|:=\bigcap_{D\sim D'\ge 0} \Supp D'
$$
is empty. So, for a free divisor $D$ the rational map $\phi_D\colon X\bir \PP^{n-1}$ associated
to $D$ in Definition \ref{def-D-morphism} is actually a morphism. The Stein factorisation of $\phi_D$
gives us a contraction $\psi_D\colon X\to Y$ such that $D\sim \psi_D^*H$ for some ample divisor $H$
on $Y$.
\end{rem}

\begin{thm}[Base point free]
Let $(X,B)$ be a klt pair where $X$ is projective. Suppose that for a nef Cartier
divisor $D$, there is some $a>0$ such that $aD-(K_X+B)$ is ample. Then,
$mD$ is free for some $m\in\N$.
\end{thm}

\begin{thm}[Rationality]
Let $(X,B)$ be a klt pair where $X$ is projective. Let $H$ be an ample Cartier divisor on $X$.
Suppose that $K_X+B$ is not nef. Then,
$$
\lambda=\max\{t>0\mid t(K_X+B)+H ~\mbox{is nef}\}
$$
is a rational number. Moreover, one can write $\lambda=\frac{a}{b}$ where $a,b\in\N$
and $b$ is bounded depending only on $(X,B)$.
\end{thm}

The proof of this is similar to the proof of the base point free theorem and
Shokurov nonvanishing theorem. So, we won't give a proof here. See [\ref{KM}] for a proof.

\begin{proof}(of Cone and Contraction theorem)
We may assume that $K_X+B$ is not nef. For any nef $\Q$-Cartier divisor $D$ define
$$
F_D=\{c\in N_1(X)\mid D\cdot c=0\}\subset \overline{NE}(X)
$$
where the inclusion follows from Kleiman ampleness criterion. Let
$$
\mathcal{C}=\overline{NE}(X)_{K_X+B\ge 0}+\overline{\sum_{D}F_D}
$$
where $D$ runs over $\Q$-Cartier nef divisors for which $\dim F_D=1$.
Suppose that $\mathcal{C}\neq \overline{NE}(X)$.
Choose a point $c\in\overline{NE}(X)$ which does not belong to $\mathcal{C}$. Now choose a rational
linear function $\alpha\colon N_1(X)\to \R$ which is positive on $\mathcal{C}-\{0\}$ but negative on
$c$. This linear function is defined by some $\Q$-Cartier divisor $G$.

If $t>0$, then $G-t(K_X+B)$ is positive on
$\overline{NE}(X)_{K_X+B=0}$ and it is positive on $\overline{NE}(X)_{K_X+B\le 0}$ for $t\gg 0$.
Now let
$$
\gamma=\min\{t>0\mid G-t(K_X+B) ~\mbox{is nef on}~ \overline{NE}(X)_{K_X+B\le 0}\}
$$

So $ G-t(K_X+B)$ is zero on some point in $\overline{NE}(X)_{K_X+B\le 0}$. Then it should be
positive on $\overline{NE}(X)_{K_X+B\ge 0}$. Therefore, $H=G-t(K_X+B)$ is ample for some rational
number $t$. Now the rationality theorem proves that
$$
\lambda=\max\{t>0\mid H+t(K_X+B) ~\mbox{is nef}\}
$$
is a rational number. Put $D=H+\lambda(K_X+B)$. Here it may happen that $\dim F_{D}>1$.
In that case, let $H_1$ be an ample divisor which is linearly independent of $K_X+B$ on $F_D$ and
let $\epsilon_1>0$ be sufficiently small. For $s>0$ let
$$
\lambda(s,\epsilon_1 H_1)=\max\{t>0\mid sD+\epsilon_1 H_1+t(K_X+B) ~\mbox{is nef}\}
$$
 
 Obviously, $\lambda(s,\epsilon_1 H_1)$ is bounded from above where the bound does not depend on $s$.
 Therefore, if $s\gg 0$, then
$$
 F_{sD+\epsilon_1 H_1+\lambda(s,\epsilon_1 H_1)(K_X+B)}\subseteq F_{D}
$$
where the inclusion is strict because $H_1$ and $K_X+B$ are linearly independent on $F_D$.
Putting all together, we get a contradiction. So, $\mathcal{C}=\overline{NE}(X)$.
Note that, by the rationality theorem, the
denominator of $\lambda(s,\epsilon_1 H_1)$ is bounded.

Now let $H_1,\dots,H_n$ be ample divisors such that together with $K_X+B$ they form a basis for
$N^1(X)$. Let $\mathcal{T}=\{c\in N_1(X)\mid (K_X+B)\cdot c=-1\}$. Let $R$ be an extremal
ray and $D$ a nef $\Q$-Cartier divisor such that $R=F_D$.
Let $c=\mathcal{T}\cap R$. Then,
$$
\lambda(D,\epsilon_j H_j)=\frac{\epsilon_j H_j\cdot c}{-(K_X+B)\cdot c}=\epsilon_j H_j\cdot c
$$
 
Therefore, this is possible only if such $c$ do not have an accumulation point in $\mathcal{T}$.
This in particular, implies that
$$
 \overline{NE}(X)=\overline{NE}(X)_{K_X+B\ge 0}+\sum_{D}R_i
$$
where $R_i$ are the $(K_X+B)$-negative extremal rays.

Now let $R$ be a $(K_X+B)$-negative extremal ray. Then, there is a nef Cartier divisor $D$
such that $R=F_D$. Thus, for $a\gg 0$, $aD-(K_X+B)$ is ample. So, by the base point free theorem,
$mD$ is free for some $m\in\N$. This gives us a contraction $\psi_{mD}\colon X\to Y$ which contracts
exactly those curves whose class belong to $R$. This proves the contractibility of $R$.
On the other hand, the morphism $mD$ is trivial only if $mD$ is ample. By construction,
$mD$ is not ample, therefore $\psi_{mD}$ is not trivial and it contracts come curve $C$.
But then the class of $C$ has to be in $R$.
\end{proof}

\begin{thm}[Kamawata-Viehweg vanishing]
Let $(X,B)$ be a klt pair where $X$ is projective. Let $N$ be an integral $\Q$-Cartier divisor on
$X$ such that $N\equiv K_X+B+M$ where $M$ is nef and big. Then,
$$
 H^i(X,N)=0
$$
for any $i>0$.
\end{thm}

For a proof of this theorem see [\ref{KM}].

\begin{thm}[Shokurov Nonvanishing]
Let $(X,B)$ be a klt pair where $X$ is projective. Let $G\ge 0$ be a Cartier divisor such
that $aD+G-(K_X+B)$ is ample for some nef Cartier divisor $D$. Then,
$$
H^0(X,mD+G)\neq 0
$$
for $m\gg 0$.
\end{thm}

\begin{rem}\label{rem-decomposition}
 Let $D$ be a Cartier divisor on a normal variety $X$. Then, we can write $D\sim M+F$
 where $M$ is movable and each component of $F$ is in $\Bs |D|$. $M$ is called
 the \emph{movable} part of $D$ and $F$ the \emph{fixed} part. There is a resolution
 $f\colon Y\to X$ such that $f^* D \sim M'+F'$ such that $M'$ is a free divisor.
\end{rem}

\begin{rem}\label{rem-ample pullback}
Every big $\Q$-divisor $D$ on a normal variety $X$ can be written as $D=H+E$
where $H$ is ample and $E\ge 0$. If in addition $D$ is also nef, then
$$
D=\frac{1}{m}((m-1)D+H)+\frac{1}{m}E=H'+\frac{1}{m}E
$$
where $m\in\N$ and $H'$ is obviously ample. So $D$ can be written as the sum of an ample $\Q$-divisor
and an effective $\Q$-divisor with sufficiently small coefficients.

 Now let $f\colon Y\to X$ be a birational contraction of normal varieties and
 $D$ an ample $\Q$-divisor on $X$. Then, there is  $E\ge 0$ contracted
 by $f$ and ample $H$ on $Y$ such that $f^*D= H+E$. Moreover, $E$ can be
 chosen with sufficiently small coefficients.
\end{rem}

\begin{rem}\label{rem-extension}
Let $S$ be a smooth prime divisor on a smooth projective variety $X$. Then, we have an exact
sequence
$$
0\to \mathcal{O}_X(-S)\to \mathcal{O}_X\to \mathcal{O}_S\to 0
$$

Moreover, if $D$ is a divisor on $X$ we get another exact sequence
$$
0\to \mathcal{O}_X(D-S)\to \mathcal{O}_X(D)\to \mathcal{O}_S(D\vert_S)\to 0
$$
which gives the following exact sequence of cohomologies
$$
0\to H^0(X, D-S)\to H^0(X, D)\toover{f} H^0(S, D\vert_S)\toover{g}
$$
$$
\hspace{2cm}  H^1(X, D-S)\toover{h} H^1(X, D)\to H^1(S, D\vert_S)
$$

It often occurs that we want to prove that $H^0(X, D)\neq 0$. The sequence above is
extremely useful in this case. However, in general $f$ is not surjective so
$H^0(S, D\vert_S)\neq 0$
does not necessarily imply $H^0(X, D)\neq 0$. But if $f$ is surjective
$H^0(S, D\vert_S)\neq 0$ implies $H^0(X, D)\neq 0$. To prove that $f$ is surjective
we should prove that  $H^1(X, D-S)=0$ or that $h$ is injective. In particular,
if $f$ is surjective and $H^0(S, D\vert_S)\neq 0$, then $S$ is not a component of the
fixed part of $D$.
\end{rem}

\begin{proof}(of the base point free theorem)
 Applying the nonvanishing theorem with $G=0$ implies that $H^0(X,mD)\neq 0$ for $m\gg 0$.
 Thus, replacing $D$ with a multiple we can assume that $H^0(X,mD)\neq 0$ for any $m\in\N$ and
 in particular, we may assume that $D\ge 0$. We will prove that $\Bs |mD|\subset \Bs |D|$
 for some $m\in \N$ where the inclusion is strict. So, applying this finitely many times
 one gets a multiple of $D$ which is a free divisor.
 
 Now suppose that $D$ is not free. By taking a log resolution
 $f\colon Y\to X$ where all the divisors involved are with simple normal crossings, we can
 further assume that $f^* D=M+F$ where $M$ is free, $F$ is the fixed part, and
 $\Supp F=\Bs |f^*D|$ by Remark \ref{rem-decomposition}. It is enough to prove that some
 component of $F$ is not in $\Bs |f^*mD|$ for some $m\in\N$.
 
 By assumptions, there is a rational number $a>0$ such that $A:=aD-(K_X+B)$ is ample. So,
 for a fixed large $m_0\in\N$ we have

$$
m_0D=K_X+B+(m_0-a)D+A
$$
and
$$
f^*((m_0-a)D+A)=E'+H
$$
where $H$ is ample and $E'\ge 0$ is exceptional over $X$ by Remark \ref{rem-ample pullback}. So,

$$
 f^*(m+m_0)D=f^*(K_X+B+mD+(m_0-a)D+A)=
$$
$$
 f^*(K_X+B+tD)+f^*(m-t)D+f^*((m_0-a)D+A)=
$$
$$
 f^*(K_X+B+tD)+f^*(m-t)D+E'+H
$$

 We can choose $t>0$ and modify $E'$ in a way that

$$
  f^*(K_X+B+tD)+E'=K_Y+B_Y+S-L
$$
where $(Y,B_Y)$ is klt, $S$ is reduced, irreducible and a component of $F$,
and $L\ge 0$ is exceptional over $X$. Therefore,

$$
f^*(m+m_0)D+L-S= K_Y+B_Y+f^*(m-t)D+H
$$

 Now by applying Kawamata-Viehweg vanishing theorem, we get

$$
 H^i(Y, f^*(m+m_0)D+L-S)=0
$$
for all $i>0$. Therefore, by Remark \ref{rem-extension}, we get the following exact sequence

$$
0\to H^0(Y, f^*(m+m_0)D+L-S)\to H^0(Y, f^*(m+m_0)D+L)\to
$$
$$
\hspace{3cm} H^0(S, (f^*(m+m_0)D+L)\vert_S)\to 0
$$

On the other hand,

$$
H^0(S, (f^*(m+m_0)D+L)\vert_S)\neq 0
$$
for $m\gg 0$ by Shokurov nonvanishing theorem because
$$
  (f^*(m+m_0)D+L)\vert_S-(K_S+B_Y\vert_S)=(f^*(m-t)D+H)\vert_S
$$
is ample. Therefore,

$$
H^0(Y, f^*(m+m_0)D+L)=H^0(X,f^*(m+m_0)D)\neq 0
$$
because $L$ is exceptional over $X$, $S$ is not a component of the fixed part of $f^*(m+m_0)D+L$
nor of $f^*(m+m_0)D$.
\end{proof}

\begin{exe}*
Prove that in the base point free theorem $mD$ is free for all $m\gg 0$.
\end{exe}

\begin{rem}[Riemann-Roch formula]
Let $D_1,\dots,D_n$ be cartier divisors on a projective variety $X$ of dimension $d$. Then,
$\mathcal{X}(\sum m_iD_i)$ is a polynomial in $m_i$ of degree $\le d$. Moreover, one can write

$$
\mathcal{X}(\sum m_iD_i)=\frac{(\sum m_iD_i)^d}{d!}+(\mbox{lower degree terms})
$$
\end{rem}

\begin{rem}[Multiplicity of linear systems]
Let $D$ be a Cartier divisor on a smooth variety $X$, and $x\in X$ a point.
Moreover, suppose that $\{h_1,\dots,h_n\}$ is a basis of $H^0(X,D)$ over $\C$. So, each
element of $h\in H^0(X,D)$ is uniquely written as

$$
h=\sum a_i h_i
$$
where we can assume that all the $h_i$ are regular at $x$ (maybe after changing $D$ in its
linear system). So, $h_i$ can be described as polynomials in $d=\dim X$ variables $t_1,\dots, t_d$
near $x$.  Thus, to give $(h)+D\ge 0$ of multiplicity $>l$ is the same as giving $h$ with
multiplicity $>l$ with respect to $t_j$. In particular, we have at most

$$
\#\{\mbox{monomials of degree $\le l$}\}=\frac{(l+1)^d}{d!}+(\mbox{lower degree terms})
$$
conditions to check.
\end{rem}

\begin{rem}\label{rem-Euler}
If two Cartier divisors $D,D'$ on a projective variety satisfy $D\equiv D'$, then
$\mathcal{X}(D)=\mathcal{X}(D')$.
\end{rem}

\begin{proof}(of Shokurov nonvanishing theorem)
We first reduce the problem to the smooth situation. Take a log resolution $f\colon Y\to X$
such that all the divisors involved have simple normal crossings. We can write

$$
f^*(aD+G-(K_X+B))=H+E'
$$
where $H$ is ample and $E'\ge 0$ is exceptional over $X$. On the other hand,

$$
 f^*(aD+G-(K_X+B))=f^*aD+f^*G-f^*(K_X+B)=
$$
$$
f^*aD+f^*G-(K_Y+B_Y-G')=f^*aD+f^*G+G'-(K_Y+B_Y)
$$
where $(Y,B_Y)$ is klt and $G'\ge 0$ is a Cartier divisor. Thus

$$
 f^*aD+f^*G+G'-(K_Y+B_Y+E')=H
$$
where $E'$ can be chosen such that $(Y,B_Y+E')$ is klt. So, we can assume that we are in the
smooth situation.

Now assume that $D\equiv 0$. Then,

$$
 h^0(X,mD+G)=\mathcal{X}(mD+G)=\mathcal{X}(G)=h^0(X,G)\neq 0
$$
by the Kawamata-Viehweg theorem and Remark \ref{rem-Euler}. So, we can assume that
$D$ is not numerically zero.

 Put $A:=aD+G-(K_X+B)$, let $k\in\N$ such that $kA$ is Cartier. Then,

$$
(nD+G-(K_X+B))^d=((n-a)D+A)^d\ge d(n-a)D\cdot A^{d-1}
$$
and since $A$ is ample and $D\neq 0$ is nef, $D\cdot A^{d-1}> 0$.

  By the Riemann-Roch theorem and the Serre vanishing theorem

$$
 h^0(X,lk(nD+G-(K_X+B)))= \mathcal{X}(lk(nD+G-(K_X+B)))= \hspace{2cm}
$$
$$
\hspace{3cm} \frac{l^dk^d(nD+G-(K_X+B))^d}{d!}+(\mbox{lower degree terms})
$$
where $d=\dim X$ and the later is a polynomial in $l$ (maybe after taking a larger $k$).
Now we can choose $n$ such that

$$
k^d(nD+G-(K_X+B))^d>2kd
$$

  On the other hand, if $x\in X-\Supp G$, for an effective divisor

$$
lkL_n\sim lk(nD+G-(K_X+B))
$$
to be of multiplicity $>2lkd$ we only need

$$
 \frac{(2lkd)^d}{d!}+(\mbox{lower degree terms of $k$})
$$

 Therefore, $L_n\sim_{\Q} nD+G-(K_X+B)$ such that $\mu_x(L_n)\ge 2\dim X$.
So, $K_X+B+L_n$ is not lc at $x$.

 Now by taking a log resolution $g\colon Y\to X$,
for $m\gg 0$ and some $t\in(0,1)$ we can write

$$
 g^*(mD+G)\equiv g^*(K_X+B+tL_n+mD+G-(K_X+B)-tL_n)\equiv
$$
$$
 g^*(K_X+B+tL_n+mD+G-(K_X+B)-tnD-tG+t(K_X+B))\equiv
$$
$$
  g^*(K_X+B+tL_n+(m-tn)D+(1-t)G+(t-1)(K_X+B)) \equiv
$$
$$
 g^*(K_X+B+tL_n+(m-tn-(1-t)a)D+(1-t)aD+(1-t)G-(1-t)(K_X+B))\equiv
$$
$$
 g^*(K_X+B+tL_n)+g^*((m-tn-(1-t)a)D+(1-t)A) \equiv
$$
$$
 g^*(K_X+B+tL_n)+H+E'\equiv
$$
$$
 K_Y+B_Y+S-G'+H+E'
$$
where $S$ is reduced and irreducible, $G'\ge 0$ is a Cartier divisor exceptional over $X$,
$H$ is ample, and $E'\ge 0$ is a suitable divisor exceptional over $X$. Hence,

$$
  g^*mD+g^*G+G'-S \equiv K_Y+B_Y+H
$$

Now by applying Kawamata-Viehweg vanishing theorem

$$
H^i(Y, g^*mD+g^*G+G'-S)=0
$$
for $i>0$. From this we deduce that

$$
 H^0(S,(g^*mD+g^*G+G')\vert_S)\neq 0 ~\implies~  H^0(Y,g^*mD+g^*G+G')\neq 0
$$

On the other hand, $g^*D\vert_S$ is a nef Cartier divisor on $S$ such that

$$
 g^*mD\vert_S+(g^*G+G')\vert_S-(K_S+B_Y\vert_S)\equiv H\vert_S
$$
is ample. Therefore, we are done by induction.
\end{proof}

\clearpage
~ \vspace{2cm}
\section{\textbf{Further studies of minimal model program}}
\vspace{1cm}

In this section we discuss aspects of minimal model program in some details.

\subsection{Finite generation, flips and log canonical models}

\begin{defn}
Let $f\colon X\to Z$ be a projective morphism of varieties and $D$ a Cartier divisor on $X$.
We call $D$ \emph{free over $Z$} if for each $P\in Z$, there is an affine open set $P\in U\subseteq Z$
such that $D$ is free on $f^{-1}U$, or equivalently if the natural
morphism $f^*f_*\mathcal{O}_X(D) \to \mathcal{O}_X(D)$ is surjective.
$D$ is called \emph{very ample over $Z$} if there is an embedding $i\colon X\to \PP_Z$ over $Z$
such that $\mathcal{O}_X(D)=i^*\mathcal{O}_{\PP_Z}(1)$. A $\Q$-Cartier divisor $D$ is called
\emph{ample over $Z$} if $mD$ is very ample for some $m\in \N$.
\end{defn}

\begin{prob}[Finite generation]\label{f.g.}
Let $f\colon X\to Z$ be a contraction of normal varieties and $D$ a $\Q$-divisor on $X$.
We associate
$$
 \mathcal{R}(X/Z,D):=\bigoplus_{m\in \Z_{\ge 0}} f_*\mathcal{O}_X(\lfloor mD\rfloor)
$$
to $D$ which is a $\mathcal{O}_Z$-algebra. When is $\mathcal{R}(X/Z,D)$ a finitely generated
$\mathcal{O}_Z$-algebra? Here by finite generation we mean that for each $P\in Z$,
there is an open affine set $P\in U\subseteq Z$ such that $\mathcal{R}(X/Z,D)(U)$ is a
finitely generated $\mathcal{O}_Z(U)$-algebra. When $Z=pt.$, we usually drop $Z$.
\end{prob}

\begin{exe}
Prove that  $\mathcal{R}(X/Z,D)$ is finitely generated iff $\mathcal{R}(X/Z,ID)$ is finitely generated
for some $I\in \N$.
\end{exe}

\begin{exe}
 Let $f\colon X'\to X$ be a contraction of normal varieties and $D'$ a divisor on $X'$ such that
 $D'=f^*D$ for a Cartier divisor $D$ on $X$. Prove that $H^0(X',mD')=H^0(X,mD)$ for
 any $m\in \N$.
\end{exe}

\begin{exe}*
Let $f\colon X\to Z$ be a contraction of normal varieties and $D$ a divisor on $X$
which is ample over $Z$. Prove that $\mathcal{R}(X/Z,D)$
is a finitely generated $\mathcal{O}_Z$-algebra. (Hint: reduce to the case $X=\PP_Z$)
\end{exe}

\begin{thm}[Zariski]
Let $f\colon X\to Z$ be a contraction of normal varieties and $D$ a Cartier divisor on
$X$ which is free over $Z$. Then, $\mathcal{R}(X/Z,D)$ is a finitely generated
$\mathcal{O}_Z$-algebra.
\end{thm}
\begin{proof}
By shrinking $Z$ we may assume that $D$ is a free divisor on $X$ and $Z$ is affine.
Now $D$ defines a contraction $\phi_D\colon X\to Z'$ such that $f$ is factored
 as $X\toover{\phi_D} Z'\to Z$ and there is a Cartier divisor $H$ on $Z'$ which is
ample over $Z$ and such that $D=f^*H$. In particular, this implies that
$$
 \bigoplus_{m\in \Z_{\ge 0}} H^0(X,mD)= \bigoplus_{m\in \Z_{\ge 0}} H^0(Z',mH)
$$
and so we are done.
\end{proof}

\begin{defn}[Log canonical ring]
Let $(X,B)$ be a pair and $f\colon X\to Z$ a contraction of normal varieties. The log canonical
ring of this pair over $Z$ is defined as

$$
 \mathcal{R}(X/Z,K_X+B)=\bigoplus_{m\in \Z_{\ge 0}} f_*\mathcal{O}_X(\lfloor m(K_X+B)\rfloor)
$$

\end{defn}

\begin{thm}
 Let $(X,B)$ be a klt pair where $X$ is projective, that is, $Z=pt.$. Then, the log canonical ring
 is finitely generated.
\end{thm}
See [\ref{BCHM}] for a proof.

\begin{thm}
 Let $(X,B)$ be a pair and $f\colon X\to Z$ a $(K_X+B)$-negative flipping contraction. Then, the flip
 of $f$ exists iff $\mathcal{R}(X/Z,K_X+B)$ is a finitely generated $\mathcal{O}_Z$-algebra.
\end{thm}
\begin{proof}
  First assume that the flip of $f$ exists: $X\to Z\leftarrow X^+$. Then,
  $\mathcal{R}(X/Z,K_X+B)=\mathcal{R}(X^+/Z,K_{X^+}+B^+)$. Since $K_{X^+}+B^+$ is ample
  over $Z$, so $m(K_{X^+}+B^+)$ is very ample over $Z$ for some $m\in\N$ and so
  $\mathcal{R}(X^+/Z,K_{X^+}+B^+)$ is a finitely generated $\mathcal{O}_Z$-algebra.

  Now assume that  $\mathcal{R}(X/Z,K_X+B)$ is finitely generated. Let $n\in\N$ be such that
  $n(K_X+B)$ is Cartier and the algebra $\mathcal{R}(X/Z,n(K_X+B))$ is locally generated by elements
  of degree one. Now define
  $X^+=\Proj \mathcal{R}(X/Z,n(K_X+B))$. Note that $f_*\mathcal{O}_X(mn(K_X+B))=\mathcal{O}_Z(mn(K_Z+B_Z))$
  where $K_Z+B_Z=f_*(K_X+B)$. By [\ref{Hart}, II-Proposition 7.10], we get a
  natural birational contraction $f^+ \colon X^+\to Z$ of normal varieties
  and a Cartier divisor $H$ on $X^+$ which is ample over $Z$. Moreover,
  $f^+_*\mathcal{O}_{X^+}(mH)=\mathcal{O}_Z(mn(K_Z+B_Z))$ for any $m\in\N$.
   It is enough to prove that $f^+$ does not contract divisors.

  Suppose that $f^+$ contracts a prime divisor $E$. Since $H$ is ample over $Z$, $E$ is not in the
  base locus of $|mH+E|$ for large $m\in\N$ after shrinking $Z$. Hence,
  $\mathcal{O}_{X^+}(mH)\subsetneq \mathcal{O}_{X^+}(mH+E)$ otherwise $E$ would be in the
  base locus of $|mH+E|$ which is not possible.
  Thus, $f^+_*\mathcal{O}_{X^+}(mH)\subsetneq f^+_*\mathcal{O}_{X^+}(mH+E)$
  for large $m\in\N$. On the other hand, since $E$ is exceptional
  $f^+_*\mathcal{O}_{X^+}(mH+E)\subseteq \mathcal{O}_Z(mn(K_Z+B_Z))$. This is a
  contradiction because $f^+_*\mathcal{O}_{X^+}(mH)=\mathcal{O}_Z(mn(K_Z+B_Z))$ by construction.
  So $f^+$ is small and, in particular, $H$ is the birational transform of $n(K_X+B)$.
\end{proof}

\begin{defn}[Log canonical model]
 Let $(X,B)$, $(Y,B_Y)$ be lc pairs and $f\colon X\bir Y$ a birational map whose inverse does not
 contract any divisors such that $B_Y=f_*B$. $(Y,B_Y)$ is called a \emph{log canonical model} for $(X,B)$
 if $K_Y+B_Y$ is ample and if $d(E,X,B)\le d(E,Y,B_Y)$ for any prime divisor $E$ on $X$ contracted by $f$.
\end{defn}

\begin{thm}
Let $(X,B)$ be a lc pair where $X$ is projective and $\kappa(X,B)=\dim X$.
Then, $(X,B)$ has a log canonical model
iff $\mathcal{R}(X,K_X+B)$ is finitely generated. Moreover, in this case the log
canonical model is given by $\Proj \mathcal{R}(X,K_X+B)$
\end{thm}
\begin{proof}(sketch)
Suppose that $(X,B)$ has a lc model $(Y,B_Y)$. Clearly,
$$
\mathcal{R}(X, K_X+B)=\mathcal{R}(Y, K_Y+Y)
$$
and since $K_Y+B_Y$ is ample, the log canonical ring is finitely generated.

Now assume that the log canonical ring $\mathcal{R}(X,K_X+B)$ is finitely generated.
Take $n\in\N$ such that $\mathcal{R}(X,n(K_X+B))$ is generated by elements of first degree.
Let $f\colon W\to X$ be a log resolution such that $f^*n(K_X+B)=M+F$ where $M$ is free
and $F$ is the fixed part. Since $M$ is free it defines a morphism $\phi_M\colon W\to \PP$.
It is well-known that the Stein factorisation of $\phi_M$ which gives the birational contraction
of $M$ is given by
$$
 W\toover{\psi_M} Y\to \phi_M(W)
$$
where $Y=\Proj \mathcal{R}(W,M)$. On the other hand,

$$
\mathcal{R}(W,M)=\mathcal{R}(X,n(K_X+B))
$$
and

$$
\Proj \mathcal{R}(X,K_X+B)=\Proj \mathcal{R}(X,n(K_X+B))
$$
so it is enough to prove that the inverse of the induced birational map $X\bir Y$
does not contact any divisors. Using more advanced tools one can prove that $F$
is contracted to $Y$ which proves the theorem.
\end{proof}

\begin{exe}\label{exe-monotonocity}
Let $(X,B)$ be a lc pair where $X$ is projective. Let $X\to Z \leftarrow X^+$ be a $(K_X+B)$-flip and
$Y$ a common log resolution of $(X,B)$ and $(X^+,B^+)$ with morphisms $g\colon Y\to X$
and $g^+\colon Y\to X^+$. Prove that for any prime divisor
$E$ on $Y$ the discrepancies satisfies

$$
 d(E,X,B)\leq d(E,X^+,B^+)
$$

Moreover, if $g(E)\subseteq \exc(f)$ then the inequality above is strict. (Hint: use the negativity lemma)
\end{exe}

\begin{exe}\label{exe-e}
Let $(X,B)$ be a klt pair. Prove that there is a log resolution $g\colon Y\to X$ such that

$$
e(X,B):=\#\{E\mid \mbox{$E$ is a prime divisor on $Y$,}
$$
$$
\hspace{4cm} \mbox{exceptional over $X$, and $d(E,X,B)\le 0$}\}
$$
is finite. Moreover, if $Y'\to X$ is another log resolution
which is over $Y$ (i.e. factors through $g$) then $e(X,B)$ is the same if it is defined using $Y'$.
\end{exe}

\begin{exe}
 Let $(X,B)$ be a klt pair. Let $X\to Z \leftarrow X^+$ be a $(K_X+B)$-flip.
 Prove that $e(X,B)\le e(X^+,B^+)$.
\end{exe}

\begin{exe}
 Let $(X,0)$ be a lc pair and $C\subseteq X$ be a subvariety of codimension $2$
 such that $C$ is not contained in the singular locus of $X$. Prove that there is
 a prime divisor $E$ on some log resolution $g\colon Y\to X$ such that $d(E,X,0)=1$.
 (Hint: use the fact that there is a log resolution $Y'\to X$ such that $Y'$ and
  $X$ are isomorphic on a nonempty open subset of $C$)
\end{exe}

\subsection{Termination of flips}

 \begin{defn}
 Let $X$ be a variety with canonical singularities and $f\colon Y\to X$ a resolution.
We define the \emph{difficulty} of $X$ as

$$
d(X):=\#\{E\mid \mbox{$E$ is a prime divisor on $Y$,}  \hspace{2cm}
$$
$$
\hspace{3cm} \mbox{exceptional over $X$, and $d(E,X,0)<1$}\}
$$
\end{defn}

\begin{exe}
Prove that $d(X)$ does not depend on the resolution. Also prove if $X\to Z\leftarrow X^+$
is a $K_X$-flip, then $d(X)\ge d(X^+)$.
\end{exe}

\begin{thm}[Termination of flips in dimension $3$]
Termination conjecture holds in dimension $3$ for varieties with terminal singularities.
\end{thm}
\begin{proof}
Let $X=X_1\bir X_2\bir X_3\bir \dots $ be a sequence of $K_X$-flips. Since $\dim X_i=3$,
in the flip diagram $X_i\to Z_i\leftarrow X_{i+1}$, each $X_i\to Z_i$ and $X_{i+1}\to Z_i$
contracts a bunch of curves. Let $C$ be a curve contracted by $X_{i+1}\to Z_i$. Since $X_{1}$
is terminal, all $X_i$ are terminal. So, by Corollary \ref{cor-terminal},
$X_{i+1}$ is smooth at the generic point of $C$, i.e., $C$ is not contained in the singular locus
of $X$. So, there is a prime divisor $E$ on a common resolution of
$X_i,X_{i+1}$ such that it is mapped to $C$ and $d(E,X_{i+1},0)=1$. Thus, $d(E,X_i,0)<1$ which in turn implies
that $d(X_i)>d(X_{i+1})$. Therefore, the sequence should terminate.
\end{proof}

 One important implication of the MMP is the following

\begin{thm}[$\Q$-factorialisation]\label{thm-Q-factorialisation}
Assume MMP in dimension $d$. Let $(X,B)$ be a klt pair of dimension $d$ where $X$ is projective.
Then, there is a small contraction $f\colon X'\to X$ such that $X'$ is $\Q$-factorial.
\end{thm}
\begin{proof}
Let $g\colon Y\to X$ be a log resolution. Define $B_Y:=B^{\sim}+E$ where $B^{\sim}$ is the birational
transform of $B$ and $E=\sum E_i$ where $E_i$ are the prime divisors contracted by $g$.
$(Y,B_Y)$ is dlt by construction. Now run the MMP over $X$, that is, in each step contract only
those extremal rays $R$ such that the class of a curve $C$ is in $R$ if $C$ is contracted to a point
on $X$. This MMP terminates and at the end we get a model $(X',B')$, a contraction $f\colon X'\to X$
for which $K_{X'}+B'$ is nef over $X$. Applying the negativity lemma proves that all the components
of $E$ are contracted in the process. This means that $f$ is small. Obviously, $X'$ is $\Q$-factorial.
\end{proof}

\begin{rem}[Adjunction]\label{rem-adjunction}
 Let $(X,B)$ be a dlt pair and $S$ a component of $\lfloor B\rfloor$. It is well-known
 that $S$ is a normal variety and we can write
 $$
  (K_X+B)\vert_S\sim_{\Q} K_S+B_S
 $$
  such that if $B=\sum b_kB_k$, then the coefficients
 of any component of $B_S$ looks like
 
 $$
 \frac{m-1}{m}+\sum \frac{l_kb_k}{m}
 $$
 for certain $m\in\N$ and $l_k\in\N\cup\{0\}$.
\end{rem}

\begin{exe}
Let $(X,B)$ be a dlt pair and $S$ a component of $\lfloor B\rfloor$. Prove that
$(S,B_S)$ is also dlt where $B_S$ is obtained by adjunction. Moreover, prove that
if $\lfloor B\rfloor$ has only one component, then $(S,B_S)$ is klt.
\end{exe}

\begin{thm}[Special termination]
Assume the MMP in dimension $d-1$. Then, any sequence of flips starting with a dlt
pair $(X,B)$ terminates near $S:=\lfloor B\rfloor$ where $X$ is projective.
\end{thm}

I will give the sketch of the proof when $S$ is irreducible. The general case
goes along the same lines.

\begin{proof}
Let $X_i\to Z_i\leftarrow X_{i+1}$ be a sequence of $(K_X+B)$-flips starting with $(X_1,B_1):=(X,B)$.
For a divisor $M$ on $X$, we denote by  $M_i$ its birational transform on $X_i$.
 Obviously, we get induced birational maps $S_i\to T_i \leftarrow S_{i+1}$ where $T_i$ is the
 normalisation of the image of $S_i$ on $Z_i$.  By adjunction
 $(K_{X_i}+B_i)\vert_{S_i}\sim_{\Q} K_{S_i}+B_{S_i}$.

We prove that $S_i\bir S_{i+1}$ and its inverse do not contract divisors for $i\gg 0$. Let $F$
be a prime divisor contracted by $S_{i+1}\to T_i$. By Remark \ref{rem-adjunction}, the coefficient
of $F$ in $B_{S_{i+1}}$ is as $b_F=\frac{m-1}{m}+\sum \frac{l_kb_k}{m}$. Now using the negativity lemma
as in Exercise \ref{exe-monotonocity}, one can prove that $d(F,S_1,B_{S_1})<-b_F$. On the other hand,
since $S$ is irreducible $(S_i,B_{S_i})$ is klt. By Exercise \ref{exe-e}, the number of such $F$
is finite. So, the inverse of $S_i\bir S_{i+1}$ does not contract divisors for $i\gg 0$.

After $\Q$-factorialising $S_i$ and lifting the sequence, we see that
$S_i\bir S_{i+1}$ does not contract divisors for $i\gg 0$ and by MMP the sequence terminates.
\end{proof}

\subsection{Minimal models}

\begin{thm}
Let $(X,B)$ be a klt pair of dimension $2$ where $X$ is
projective. Then, $(X,B)$ has a unique minimal model, or a Mori fibre space (not unique).
\end{thm}
\begin{proof}
If $K_X+B$ is not nef, then by the cone and contraction theorem, there is a
$(K_X+B)$-negative extremal ray $R$. Let $f\colon X\to Z$ be the contraction of
$R$. If $R$ is of fibre type, then we get a Mori fibre space and we are done.
Otherwise $f$ is birational. Since $X$ is a surface, $f$ cannot be a flipping
contraction, so $f$ is a divisorial contraction which is nontrivial. It is easy to see that
$\rho(X)>\rho(Z)$. Now replace $(X,B)$ by $(Z,f_*B)$ and continue the
process. Obviously, the program stops after finitely many steps. So, we end up
with a minimal model or a Mori fibre space.

The uniqueness of the minimal model follows from the more general
Theorem \ref{thm-minimal models}.
\end{proof}

In the proof of the previous theorem we did not explain why
$K_Z+f_*B$ is $\Q$-Cartier. Theorem \ref{thm-Q-factorialisation} shows that $Z$ is actually
$\Q$-factorial because we can first prove MMP in dimension $2$ for $\Q$-factorial pairs
and use the proof of Theorem \ref{thm-Q-factorialisation} to show that $Z$ is in general
$\Q$-factorial.

\begin{exe}
 Let $(X,B)$ be a klt pair of dimension $2$ where $X$ is
projective and $\kappa(X,B)\ge 0$. Prove that $(X,B)$ has a minimal
model.
\end{exe}

\begin{exe}
 Let $X$ be a smooth projective surface. Prove that the minimal or Mori fibre spaces
 obtained in the previous theorem are also smooth if $B=0$.
\end{exe}

\begin{thm}
Let $(X,B)$ be a klt pair of $\kappa(X,B)=\dim X$ which is minimal, that is, $K_X+B$ is nef.
Then, the abundance conjecture holds in this case.
\end{thm}
\begin{proof}
Since $\kappa(X,B)=\dim X$, $K_X+B$ is big. It is also nef by assumption.
So, we can write $K_X+B\sim_{\Q}H+E$ where $H$ is ample and $E\ge 0$ has sufficiently
small coefficients. In particular,
$r(K_X+B)-(K_X+B+E)=(r-2)(K_X+B)+H$ is ample if $r\ge 2$. So, applying the base
point free theorem proves that $m(K_X+B)$ is base point free for some $m\in\N$.
\end{proof}
 
\begin{defn}[Flop]
Let $(X,B)$ be a pair where $X$ is projective. A flop $X\to Z\leftarrow X^+$ is just like a flip except that
we assume that $(K_X+B)\equiv 0$ and $K_{X^+}+B^+\equiv 0$ over $Z$.
\end{defn}

\begin{thm}\label{thm-minimal models}
Let $(X,B)$ be a klt pair where $X$ is projective. Let $(Y_1,B_1), (Y_2,B_2)$ be two minimal
models of $(X,B)$. Then, $Y_1$ and $Y_2$ are isomorphic in codimension $1$. Moreover,
assuming the MMP, the induced birational map $Y_1\bir Y_1$ is decomposed into
 finitely many flops.
\end{thm}
\begin{proof}
Since $Y_1,Y_2$ are minimal models of $(X,B)$, $K_{Y_1}+B_1$ and $K_{Y_2}+B_2$ are nef.
Now take a common log resolution $f_1\colon Y\to Y_1$, $f_2\colon Y\to Y_2$. Apply the negativity
lemma to prove that
$$
f_1^*(K_{Y_1}+B_1)=f_2^*(K_{Y_2}+B_2)
$$
which implies that $Y_1\bir Y_2$ is an isomorphism in codimension one.

Let $H_2\ge 0$ be an ample divisor on $Y_2$ and $H_1$ its birational transform on $Y_1$.
Take $t>0$ such that $K_{Y_1}+B_1+tH_1$ is klt. If $0<s<t$ is small enough, we can find
an extremal ray $R$ such that $(K_{Y_1}+B_1+sH_1)\cdot R<0$ and $(K_{Y_1}+B_1)\cdot R=0$
unless $Y_1\bir Y_2$ is an isomorphism.
In particular, $(K_{Y_1}+B_1+tH_1)\cdot R<0$. Since $K_{Y_1}+B_1$ is nef, the contraction
of $R$ is flipping. Now by MMP the flip $Y_1\to Z_1\leftarrow Y_1^+$ for $R$ exists and
$(Y_1^+,B_1^+)$ is also a minimal model of $(X,B)$. By continuing  this process we get a
sequence of  $(K_{Y_1}+B_1+sH_1)$-flips which terminates by MMP. In other words, we arrive
at a model where the birational transform of $K_{Y_1}+B_1+sH_1$ is nef. Now since
$K_{Y_2}+B_2+sH_2$ is ample, we should arrive at $Y_2$.
\end{proof}

\subsection{Mori fibre spaces}

\begin{thm}
Let $(X,B)$ be a pair which has a Mori fibre space where $X$ is projective.
Then, $\kappa(X,B)=-\infty$.
\end{thm}
\begin{proof}
Le $(Y,B_Y)$ be a Mori fibre space for $(X,B)$ and $f\colon Y\to Z$ the corresponding
contraction which is the contraction of a $(K_X+B)$-negative extremal ray. Suppose
that $\kappa(X,B)\ge 0$, i.e., there is $M\ge 0$ such that $K_X+B\sim_{\Q}M$. If
$M_Y$ is the birational transform of $M$ on $Y$, then $K_Y+B_Y\sim_{\Q}M_Y\ge 0$.
So, $M_Y\cdot C<0$ for any curve $C$ contracted by $f$. This is not possible
because if we take $C$ such that $M_Y$ does not contain it, then $M_Y\cdot C\ge 0$,
a contradiction. Therefore, $\kappa(X,B)=\infty$.
\end{proof}

Note that the minimal model conjecture claims that the inverse is also true, that is,
if $\kappa(X,B)=-\infty$, then $(X,B)$ has a Mori fibre space.

\begin{thm}
Let $(X,B)$ be a pair where $X$ is a projective surface and $f\colon X\to Z$ a Mori fibre space
with $\dim Z=1$. Then, a general fibre of $f$ is isomorphic to $\PP^1$. Moreover, if $X$ is
smooth, then every fibre is isomorphic to $\PP^1$.
\end{thm}
\begin{proof}
Since $X$ is normal, it has only finitely many singular points, so almost all the
fibres are in the smooth part of $X$. Now let $F=f^*P$ and $G=f^*Q$ where $P,Q\in Z$
are distinct and $F$ is in the smooth locus of $X$. Let $C$ be a component of $F$.
Then, $G\cdot C=0$. If $F$ has more than one component, then there is another component $C'$
which intersects $C$ because $f$ is a contraction and it has connected fibres.
So, $C\cdot C'>0$. This is impossible since $f$ is an extremal contraction and the class
of all the curves contracted by $f$ are in the same ray. Therefore, $F$ has only one component.
In particular, it means that $C^2=0$.

So, we can assume that $C$ is the single component of $F$. Thus we have

$$
(K_X+B+C)\cdot C=(K_X+C)\cdot C+B\cdot C=2p_a(C)-2+B\cdot C<0
$$
which implies that $2p_a(C)-2<0$ and the arithmetic genus $p_a(C)=0$ if $B\cdot C\ge 0$.
Of course, $B\cdot C\ge 0$ holds for almost all fibres $F$.
In this case, since $C$ is irreducible, $C\simeq \PP^1$.
\end{proof}

\begin{defn}[Fano pairs]
A lc pair $(X,B)$ is called a Fano pair if $-(K_X+B)$ is ample.
\end{defn}

\begin{exe}
Let $(X,B)$ be a Fano pair of dimension $1$. Prove that $X\simeq \PP^1$ and
$\deg B <2$.
\end{exe}

\begin{rem}
Let $(X,B)$ be dlt pair and $f\colon X\to Z$ be a Mori fibre space. Then,
$(F,B_F)$ is a Fano pair for a general fibre $F$ where $K_F+B_F\sim_{\Q} (K_X+B)\vert_F$.
\end{rem}

\begin{rem}[Sarkisov program]
Let $(X,B)$ be a lc pair and $(Y_1,B_1)$ and $(Y_2,B_2)$ two Mori fibre spaces
for $(X,B)$ with corresponding contractions $f_1\colon Y_1\to Z_1$ and
$f_2\colon Y_2\to Z_2$. Then, it is expected that the induced birational map
$Y_1\bir Y_2$ is decomposed into simpler ones which are called Sarkisov links.
For example, some links look like

$$
\xymatrix{ &Y'\ar[ld]^f\ar[rd]^g&\\
Y_1\ar[d]& \dashrightarrow & Y_2\ar[d] \\
Z_1&&Z_2}
$$
where $f,g$ are divisorial contractions.
\end{rem}

\begin{exa}
Let $X=\PP^1\times\PP^1$ and $f_1\colon X\to \PP^1$ and $f_2\colon X\to \PP^1$ the two
natural projections. Both $f_1$ and $f_2$ are Mori fibre spaces and for any fibre $F$
of $f_1$ or $f_2$, $K_X\cdot F=-2$. So, even on a single variety there could be different
Mori fibre structures. Now let $P\in X$ and $g\colon Y\to X$ be the blow up of $X$ at $P$.
If $F_i$ is the fibre of $f_i$ passing through $P$, then $F_i^2=0$ but $\tilde{F}_i^2=-1$
on $Y$ where $\sim$ denotes the birational transform. Thus, $\tilde{F}_1$ and  $\tilde{F}_2$
are disjoint $-1$-curves. Let $h\colon Y\to X'$ be the contraction of both $\tilde{F}_1$
and  $\tilde{F}_2$. It turns out that $X'\simeq \PP^2$ and this shows that $Y$ has at least three
different Mori fibre spaces.
\end{exa}

\begin{exe}
A smooth projective surface $X$ is called \emph{ruled} if it is birational to $\PP^1\times C$
for some curve $C$. Prove that any ruled surface has a Mori fibre space. In particular,
$\kappa(X)=-\infty$ for a ruled surface.
\end{exe}

It is well-known that if $X$ is a smooth projective surface and $\kappa(X)=-\infty$, then
$X$ is ruled.

\clearpage
~ \vspace{2cm}
\section{\textbf{Appendix}}
\vspace{1cm}

\begin{enumerate}
\item \emph{Riemann-Roch for curves;} $\mathcal{X}(\mathcal{O}_X(D))=1+\deg D-p_a(X)$.\\
\item \emph{Riemann-Roch for surfaces;} $\mathcal{X}(\mathcal{O}_X(D))=\frac{1}{2}D\cdot (D-K)+1+p_a(X)$.\\
\item \emph{Adjunction on surfaces;} For an effective divisor $D$ on a smooth projective surface $X$
we have $p_a(D)-1=\frac{1}{2}D\cdot (K_X+D)$.\\
\end{enumerate}

\vspace{1cm}

\end{document}